\documentclass[a4paper,english,12pt]{article}
\usepackage[logonly]{trace}
\usepackage{babel}
\usepackage{amsfonts, amsmath, amssymb}
\usepackage{graphicx}
\usepackage{pstricks}
\usepackage{psfig}
\usepackage{multirow}
\usepackage{fancyhdr}
\usepackage{vmargin, fancybox}

\newcommand{\mathsym}[1]{{}}

\def\N{\textrm{I\kern-0.21emN}}

\def\R{\textrm{I\kern-0.21emR}}
\def\Q{\textrm{l\kern-0.5emQ}}

\begin{document}

\title{A note on the Compound Burgers-Korteweg-de
Vries Equation with higher-order nonlinearities and its traveling
solitary waves}

\author{Claire David \footnotemark[1] \\
 \\
\footnotemark[1] Universit\'e Pierre et Marie Curie-Paris 6  \\
Institut Jean Le Rond d'Alembert, UMR CNRS 7190 \\
Bo\^ite courrier $n^0162$, 4 place Jussieu, F-75252 Paris cedex
05, France\\}

\date{}

\maketitle

\begin{abstract}
In this paper, we study a compound Korteweg-de Vries-Burgers
equation with a higher-order nonlinearity.
  A class of solitary wave solutions is
obtained by means of a series expansion.
\end{abstract}

\section{Introduction}
\label{sec:intro}

\noindent Consider either the a Korteweg-de Vries-Burgers-type
equation of the following form:

\begin{equation}
\label{CBKDV}
u_t + \alpha u^p\, u_x + \beta\,u^{2\,p}\,u_x+\gamma \,
u_{xx}+\mu \, u_{xxx} = 0 \end{equation}

\noindent where $\alpha$, $\beta$, $\gamma$, $\mu$ and $s$ are
real constants, while $p$ is a positive number, or, more
generally:

\begin{equation}
u_t + P(u)\, u_x + \beta\,u^{2\,p}\,u_x+\gamma \, u_{xx}+\mu \,
u_{xxx} = 0 \end{equation}

\noindent where $P$ is a generalized polymon, i.e. of the form:

\begin{equation}
P(u)=\displaystyle \sum_{r\in {\R}^{+}} \alpha_r\,u^r
\end{equation}

\section{Traveling solitary wave solutions}

\noindent Assume that equation (\ref{CBKDV}) has the solution of
the form:

\begin{equation}
\label{sol} u(x,t)=u(\zeta) \,\, \, , \, \, \, \zeta=x-v\,t
\end{equation}

\noindent where $v$ is the velocity.\\
\noindent Substituting it into (\ref{CBKDV}) yields:

\begin{equation}
-v\,u'(\zeta) + \alpha u^p\, u'(\zeta) +
\beta\,u^{2\,p}\,u'(\zeta)+\gamma \, u''(\zeta)+\mu \,
u^{(3)}(\zeta) = 0
\end{equation}

\noindent Performing one integration, we have then:

\begin{equation}
\label{CBKDV_int} -v\,u(\zeta) + \frac {\alpha}{p+1}
{u(\zeta)}^{p+1} + \frac {\beta}{2p+1} {u(\zeta)}^{2p+1}+\gamma \,
u'(\zeta)+\mu \, u^{(2)}(\zeta)+ d = 0
\end{equation}
\noindent where $d$ is an integration constant.\\
\noindent For sake of simplicity, we shall take $d=0$.

\noindent Set:

\begin{equation}
a=\frac {\alpha}{\mu(p+1)} \,\,\, , \,\,\, b=\frac
{\beta}{\mu(2p+1)} \,\,\, , \,\,\,  c = \frac {\gamma} {\mu} \,\,\,
, \,\,\,  r=\frac {v} {\mu}
\end{equation}

\noindent Equation (\ref{CBKDV_int}) can thus be written as:

\begin{equation}
\label{CBKDV_int2} -r\,u(\zeta) + a\, {u(\zeta)}^{p+1} + b\,
{u(\zeta)}^{2p+1}+c \, u'(\zeta)+ u^{(2)}(\zeta) = 0
\end{equation}

\noindent Consider the surface $S_u$ in the three-dimensional
euclidean space:

\begin{equation}
\label{Surf} -r\,X + a\, X^{p+1} + b\, X^{2p+1}+c \, Y+ Z = 0
\end{equation}

\noindent $u$, $u'$, $u^{(2)}$ are traced on this surface. The
knowledge of a parametrization of this surface will thus lead to
the determination of $u$.

\subsection{The integer case}

\noindent Contrary to previous works (\cite{David}), there is no
useful information available about the solutions or their
profiles. Thus, we choose to search the solution as the following
series expansion:
\begin{equation}
\label{Serie} u(x, t) = \sum_{k = 0}^{+\infty} U_k\,e^{k\,z}
\end{equation}
where the $U_k's$ are constants to be determined.
\\

\noindent Substitution of (\ref{Serie}) into equation
(\ref{CBKDV}) leads to an equation of the following form:
\begin {equation}
\sum_{k=0}^{+\infty} P_k(U_k,a,b,c,r) \,e^{k\,z}=  0,
\end{equation}
\noindent where the $P_k$ $(k=0,\,...,\,+\infty)$, are polynomial
functions of the $U_k$ and of $a$, $b$, $c$, $r$. The solution is
obtained equating the $P_k$ $(k=0,\,...,\,,+\infty)$ to zero.\\

\noindent With the aid of mathematical softwares such as
Mathematica, the previous system can be solved.\\

\noindent In the following, we present the solution obtained in
the case $p=1$, for $r = 1$, $a = 0.4$, $b = 0.01$, $c = 0.2$. The
series (\ref{Serie}) is truncated at $n=3$:

\begin{equation}
 u(x, t) = \sum_{k = 0}^{3} U_k\,e^{k\,z}
\end{equation}

\noindent The values of the coefficients are:

$$\left \lbrace \begin{array}{ccc}% autant de c que de colonnes
B_0 &=& -0.2523887531009444 \\ B1 &=&
 7.920512040580792 +
  16.296799786819456\,i  \\
B_2 &=&
 24.87642134042838 - 31.6589105912486\,i \\
  B_3 &=& -59.69562063336409 - 12.94860377480183\,i
\end{array}
 \right .$$

\begin{figure}[h!]
\centerline{
 \includegraphics[width=8cm,height=7cm]{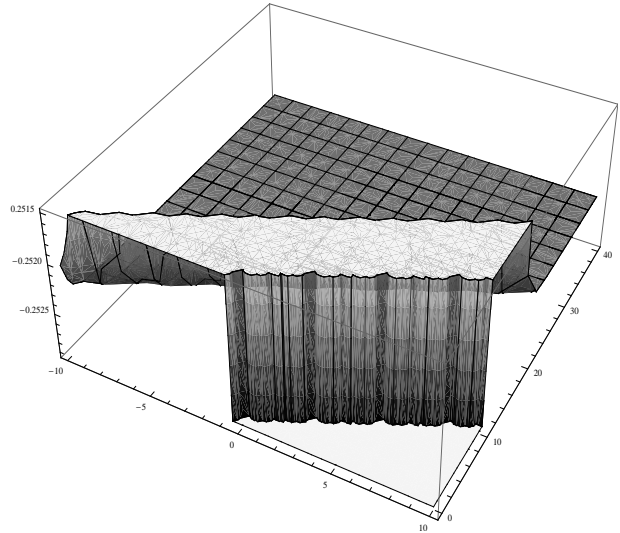}}
  \caption{\small{The real part of the solitary wave solution in the case surface in the case $p=1$, for $r = 1$, $a = 0.4$, $b = 0.01$, $c = 0.2$.}} \label{Sol}
\end{figure}

\noindent Figure \ref{Sol} displays the real part of the solitary
wave solution.

\subsection{The non-integer case}

\noindent In the specific case where the real number $p$ is not an
integer, the exact determination of the traveling solitary wave
solution of (\ref{CBKDV}) becomes impossible. Yet, by means of
numerical methods of surfaces reconstruction, one can approximate
$u$. Also, plotting the surface $S_u$ can yield interesting
informations on $u$.

\begin{figure}[h!]
\centerline{
 \includegraphics[width=7cm,height=5cm]{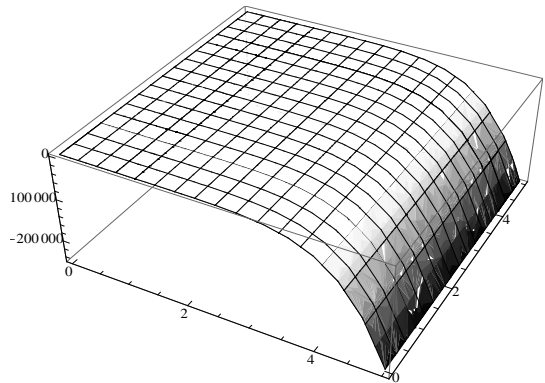}}
  \caption{\small{The surface in the case $a=1$, $b=2$, $c=3$, $r=1$, $p=\pi$.}} \label{SurfPi}
\end{figure}

\noindent \addcontentsline{toc}{section}{\numberline{}References}

\end{document}